\documentstyle[12pt]{article}
\newcommand{\bk}{\\[0.01in] \hspace*{0.5in}}
\newcommand{\Rm}{Riemannian metric }
\newcommand{\Rms}{Riemannian metrics }
\newcommand{\qed}{ \ \ \ \ \ {\bf Q.E.D.}\\[0.1in]}
\newcommand{\og}{{\overline g}}
\newcommand{\oR}{{\overline R}}
\newcommand{\oM}{{\overline M}}
\begin{document}

\title{Conformal deformation of warped products and scalar curvature 
functions on
open manifolds}
\author{Man Chun LEUNG\\ Department of Mathematics,\\ National University of
Singapore,\\ Singapore 119260\\ matlmc@nus.sg}
\date{April 1996}
\maketitle
\begin{abstract}We discuss conformal deformation and warped products 
on some open manifolds. We discuss how these can be applied to construct
\Rms with specific scalar curvature functions.  
\end{abstract}
\vspace{0.2in} KEY WORDS:  scalar curvature, warped product, conformal
deformation
\\[0.075in]   1991 AMS MS Classifications: 53C20, 53A30

\vspace{0.5in}

{\bf \Large 1. \ \ Introduction}

\vspace{0.3in}

In this paper we study scalar curvature functions on some open manifolds. A
classification result of Kazdan and Warner (with an improvement by B\'erard
Bergery)  states that if
$N$ is a compact
$n$-manifold without boundary,
$n \ge 3$, then $N$ belongs to one of the following three catagories ([2] p.
125).\\[0.05in] {\bf (A)} \ \ Any smooth function on $N$ is the scalar 
curvature
of some Riemannian metric on $N$.\\[0.01in]
{\bf (B)} \ \ A smooth function on $N$ is the scalar curvature of some
Riemannian metric on $N$ if and only if the function is either 
identically zero
or strictly negative somewhere; moreover, any metric with vanishing scalar
curvature is Ricci-flat.\\[0.01in]
{\bf (C)} \ \ A smooth function on $N$ is the scalar curvature of some
Riemannian metric on $N$ if and only if the function is negative
somewhere.\bk
Thus if the manifold $N$ admits a Riemannian metric of
positive scalar curvature, then it belongs to class (A). If $N$ admits a
Riemannian metric of zero scalar curvature and cannot admit a metric of 
positive
scalar curvature, then it belongs to class (B). And any Riemannian metric 
on a
manifold belongs to class (C) has scalar curvature negative somewhere. This
completely answers the question of which smooth functions  are scalar 
curvatures
of Riemannian metrics on a compact manifold $N$. For noncompact 
manifolds, many
important works have been done on the question of how to determine which 
smooth
functions are scalar curvatures of complete \Rms on an open manifold. 
Results of
Gromov and Lawson [5] show that some open manifolds cannot carry complete
Riemannian metrics of positive scalar curvature, for examples, weakly 
enlargeable
manifolds. Furthermore, they show that some open manifolds cannot even admit
complete Riemannian metrics with scalar curvatures uniformly positive 
outside a
compact set and with Ricci curvatures bounded ([5], [6] p. 322). On the other
hand, it is known that each open manifold of dimension bigger than 2 
admits a
complete Riemannian metric of constant negative scalar curvature [3]. It
follows from the results of Aviles and McOwen [1] that any bounded negative
function on an open manifold of dimension bigger than 2 is the scalar 
curvature
of a complete Riemannian metric.\bk  
By using conformal deformation, Ni and other authors have studied
which functions on ${\bf R}^n$ are scalar curvatures of complete Riemannian
metrics that are conformal to the Euclidean metric [8]. While in [10], Ratto,
Rigoli and V\'eron give a rather detailed study of similar problem on a
hyperbolic space using conformal deformation. In [7], we apply conformal
deformation to study scalar curvatures on more general types of open 
manifolds.
Due to large  variety of structures on open manifolds, it is rather 
unclear how
to  consider the scalar curvature question on an arbitrary open manifold. We
will mainly restrict ourselves to open manifolds that have 
compactifications.\bk
Let ${\overline M}$ be a compact $(n+1)$-manifold with boundary $\partial 
M$ and
interior $M$. The boundary $\partial M$ has finite number of connected
components, and each component is a compact $n$-manifold without 
boundary. In
this paper, we discuss the method of using warped products and conformal
deformation to construct complete Riemannian metrics on
$M$ with specific scalar curvatures. By making use of the boundary, we can
construct warped products at the ends of
$M$. It is shown that if the boundary components belong
to class (A) or (B), then $M$ admits a complete \Rm with positive scalar
curvature outside a compact set. If a boundary component belongs to class 
(C),
then we can construct complete \Rm with scalar curvature approaching zero 
near
the end. We discuss restrictions of using warped products and conformal
deformation to obtain complete \Rm of positive scalar curvature outside a
compact set. If a connected component of
$\partial M$ belongs to class (C), then we show that, under mild 
assumptions of
warping functions, it is not possible to conformally deform warped product
metrics to complete metrics on $M$ with nonnegative scalar curvature 
outside a
compact set (theorem 3.8). In section 4 we discuss the scalar curvature
functions and conformal deformation of more general type of \Rms known as 
polar
type \Rms. We discuss some restrictions of using polar type \Rms and
conformal deformations to obtain nonnegative scalar curvature outside a 
compact
set of $M$ (theorem 4.18).

\vspace{0.3in}

{\bf \Large 2. \ \ Boundary components in class (A) or (B)}

\vspace{0.3in}

Let $(N, g)$ be a Riemannian manifold of dimenison $n$ and let
$f: (2, \infty) \to {\bf R}^+$ be a smooth function. The warped product 
of $N$
and
$(2, \infty)$ with warping function $f$ is defined to be the Riemannian 
manifold
$((2, \infty) 
\times N, g')$ with
$$g' = dt^2 + f^2 (t) g\,. \leqno (1.1)$$ Let $R (g)$ be the
scalar curvature of $(N, g)$. Then the scalar
curvature
$R$ of $g'$ is given by the equation
$$
R (t, x) = {1\over {f^2 (t)}} \{ R (g) (x) - 2nf (t) f'' (t) - n (n - 1)
|f' (t)|^2\}
\leqno (1.2)
$$
for $t \in (2, \infty)$ and  $x \in N\,.$ In section 4 we discuss a more
general formular for the scalar curvature of the
\Rm $dt^2 + f^2 (x, t) g$, where $f$ is a positive smooth function of $(2,
\infty)$ and
$N$. When $f$ is a constant function, the metric $g'$ in (1.1) is known
as a conic metric. See [7] for a discussion of scalar curvature and conformal
deformation of conic metrics. If we denote  
$$u (t) = f^{{n + 1}\over 2} (t)\,, \ \ \ \ t > 2\,,$$  
then equation (1.2) can be changed into [4]
$${{4n}\over {n + 1}} u'' + R u - R (g) u^{{n - 3}\over {n + 1}} =  0 \,.
\leqno (1.3)$$
Let ${\overline M}^n$ be a compact $(n + 1)$-manifold with boundary 
$\partial M$
and interior $M$. In this paper we assume that the boundary is nonempty 
and has
connected components $N_1, N_2,..., N_k$. Each $N_i$ is a compact 
$n$-manifold
without boundary. A neighborhood of $N_i$ in $M$ is diffeomorphic to $(2,
\infty) \times N_i$, $i= 1,2,..,k$. The existence of a complete \Rm on $M$
with constant negative scalar curvature has been proved in [3]. In the 
present
case we have a simplier proof.
\\[0.1in]
{\bf Proposition 1.4.} \ \ {\it If $n \ge 3$, then on $M$ there is a complete
metric of constant negative scalar curvature which is a product
metric near infinity.}\\[0.1in]  {\bf Proof.} \ \ On the compact manifold
$N_i$, there exists a metric $g_i$ of constant negative scalar curvature.
Using the product metric on each $(2, \infty) \times N_i$ and extending to
the whole $M$, we obtain a complete \Rm $g_1$ on $M$ with scalar
curvature $R (g_1)$ being a  negative constant outside a compact set. Let 
$u$ be
a nonnegative smooth function on
$M$ such that $u \equiv 1$ on $\partial M \times (2, b)$ for some $b > 
2$, $u
\equiv 0$ on $\cup_{i = 1}^k (b+ 1, \infty) \times N_i$ and
$|\bigtriangledown u| < C_o$ for some positive constant $C_o$ independent on
$b$. For $b$ large enough, we have
$$\int_M ( {{4 n}\over {(n - 1)}} |\bigtriangledown u|^2 + R (g_1) u^2) 
dv_g <
0\,.$$  
Using a result of Aviles and McOwen [1], $g$ can be conformally
deformed into a complete Riemannian metric of constant negative scalar
curvature. Furthermore, the conformal factor can be chosen to be equal to a
positive constant outside a compact set.\qed 
{\bf Corollary 1.5. [1]} \ \ {\it Any
negative smooth function on $\overline M$ is the scalar curvature of a 
complete
metric on $M$.}\\[0.1in]
{\bf Proof.} \ \ The function is a bounded negative function. Using
the complete \Rm constructed above, we can conformally deform it into a
complete \Rm with the prescribed scalar curvature.\qed
\hspace*{0.5in}Consider the case when scalar curvature functions can be 
positive. For
$1
\le i
\le k,$ if the manifold
$N_i$ admits a
\Rm  of positive scalar cuvature, then the product metric
$$dt^2 + g_i$$  
has positive scalar curvature on $(2\,, \infty) \times N_i$. If $N_i$ 
admits a
\Rm of zero scalar curvature, then we let $u (t) = t^\alpha$ in (1.3), 
where 
 $\alpha \in (0, 1)$ is a constant. We have 
$$R = {{4n}\over {(n + 1)}} \alpha (1 - \alpha) 
{1\over {t^2}} > 0\,, \ \ \ \ t > 2\,.
\leqno (1.6)$$
Therefore we obtain the following.\\[0.1in]
{\bf Theorem 1.7.} \ \ {\it For $n \ge 3$, let $M$ be the interior of a 
compact
$(n + 1)$-manifold with boundary. Suppose that the boundary components 
are in
class (A) or (B), then on $M$ there is a complete \Rm 
of positive scalar curvature outside a compact set}\\[0.1in]
\hspace*{0.5in}
We note that the term $\alpha (1 - \alpha)$ achieves its maximum when
$\alpha = 1/2$. And when $u = t^{1\over 2}$ we have 
$$R = {{4n}\over {(n + 1)}} {1 \over 4} {1\over {t^2}}\,, \ \ \ \ t > 2\,.$$
We show that this is almost the best possible.\\[0.1in]
{\bf Lemma 1.8.} \ \ {\it If $R (g) = 0$, then there are no
positive solutions to the equation (1.3) with}
$$R (t) \ge {{4n}\over {(n + 1)}} {c \over 4} {1 \over {t^2}} \ \ \ \ 
{for} 
\ \ \ t \ge t_o\,,$$   
{\it where $c > 1$ and $t_o > 2$ are constants.}\\[0.1in]
{\bf Proof.} \ \ Assume that 
$$R (t) \ge {{4n}\over {(n + 1)}} {c \over 4} {1 \over {t^2}} \ \ \ \ 
{\rm for} 
\ \ \ t \ge t_o\,,$$   
with $c > 1$. Equations (1.3) gives
$$t^2 u'' (t) + {c\over 4} u \le 0\,.$$ 
Let 
$$u (t) = t^\alpha v (t)\,, \ \ \ t \ge t_o\,,$$ 
where $\alpha > 0$ is a constant and $v (t) > 0$ is a smooth function. 
Then we
have
$$u'' (t) = \alpha (\alpha - 1) t^{\alpha - 2} v (t) + 2 \alpha t^{\alpha 
- 1}
v' (t) + t^\alpha v'' (t)\,.$$ 
And we obtain
$$t^\alpha v (t) [\alpha (\alpha - 1) + {c\over 4}] + 2 \alpha t^{\alpha 
+ 1} v'
(t) + t^{\alpha + 2} v'' (t) \le 0\,. \leqno (1.9)$$ 
Let $\delta$ be a positive
constant such that $\delta^2 = (c - 1)/4$. Then we have
$$\alpha (\alpha - 1) + {c\over 4} = (\alpha - {1\over 2})^2 + {{c - 
1}\over 4}
\ge \delta^2\,.$$
Then $\delta$ is a constant independent on $\alpha$. (1.9) gives
$$2 \alpha t v' (t) + t^2 v'' (t) \le - \delta^2 v (t)\,. \leqno (1.10)$$ 
Let $\beta = 2
\alpha$ and we choose $\alpha > 0$ such that $\beta < 1$, that is,
$\alpha < 1/2$. Then (1.10) becomes
$$(t^\beta v' (t))' \le -{{\delta^2 v (t)}\over {t^{2 - \beta} }}\,.$$ 
Upon integration we have
$$t^\beta v' (t) - \tau^\beta v' (\tau) \le - \int_\tau^t {{\delta^2 v 
(s)}\over
{s^{2 - \beta} }} ds\,, \ \ \ \ \ t > \tau > t_o\,. \leqno (1.11)$$ 
If $v' (\tau) \le 0$ for some $\tau > t_o$, then (1.11) impliest that 
$$t^\beta v' (t) \le - C$$ 
for some positive constant $C$. We have
$$v (t) \le v (\tau) - \int_\tau^t {C\over {s^\beta}} ds = v (\tau) - 
{{t^{1 -
\beta} }\over { 1 - \beta}} \vert^t_\tau \to - \infty\,,$$ as $\beta < 
1$. Hence
$v (t) < 0$ for some $t$, contradicting that $u (t) > 0$ for all $t \ge t_o$.
Thus we have $v' (t) > 0$ for all $t > t_o$. (1.11) implies that
$$\tau^\beta v' (\tau) - \int_\tau^t {{\delta^2 v (s)}\over {s^{2 - 
\beta} }} ds
\ge 0$$ 
for all $t > \tau > t_o$. As $v' (t) > 0$ for all $t > t_o$, we have
$$\tau^\beta v' (\tau) \ge v (\tau)  \int_\tau^t {{\delta^2 }\over {s^{2 -
\beta} }} ds = v (\tau) {1\over {s^{1 - \beta} }} [ - {{\delta^2}\over { 
1 -
\beta}} ] |^t_\tau\,.$$ 
Let $t \to \infty$ we have
$$\tau^\beta v' (\tau) \ge {{v (\tau)}\over {\tau^{1 - \beta} }} 
{{\delta^2}\over
{ 1 - \beta}}\,.$$ 
Or after changing the parameter we have 
$${{v' (t)}\over {v (t)}} \ge {1\over t} {{\delta^2}\over { 1 - \beta}}\,,
\ \ \ t > t_o\,.$$ 
Choosing $\alpha < 1/2$ close to $1/2$ so that $\beta < 1$ is
close to
$1$ and using the fact that $\delta$ is independent on $\alpha$ or 
$\beta$, we
have
$${{v' (t)}\over {v (t)}} \ge {N\over t}$$ for a big integer $N > 2$. 
This gives
$$v (t) \ge C t^N\,, \ \ \ \ t > t_o\,,$$ 
where $C$ is a positive constant. (1.11) implies that
$$t^\beta v' (t) \le \tau^\beta v' (\tau)  - \int_\tau^t {{C \delta^2 
s^N}\over
{s^{2 - \beta} }} ds \to -\infty \ \ \ \ {\mbox{as}} \ \ \ t \to
\infty\,.$$ Thus $v' (t) < 0$ for $t$ large.\qed 
\hspace*{0.5in}In particular, if $R (g) = 0$, then using warped product 
it is
impossible to obtain a \Rm of uniformly positive scalar curvature. The best
we can do is when $u (t) = t^{1\over 2}$, or $f (t) = t^{1\over {n + 1}}$,
where the scalar curvature is positive but goes to zero at infinity. It 
can also
be shown that one cannot conformally deform the metric 
$$dt^2 + t^{2\over {n + 1}} g$$ into a complete metric of uniformly positive
scalar curvature outside a compact set. For simplicity we assume that the
boundary of
$\overline M$ is connected. The case where $\partial \overline M$ has 
more than
one connected component is similar.\\[0.1in]
{\bf Theorem 1.12} \ \ {\it For $n \ge 3$, let $\overline M$ be a compact
$(n + 1)$-manifold with boundary $N$. Suppose that $N$ admits a \Rm $g$ 
of zero
scalar curvature. Let $g'$ be a complete \Rm on $M$ with} 
$$g' (t, x) = dt^2 + t^{2\over {n + 1}} g (x) \ \ \ \ {\mbox {on}} \ \ \ (2,
\infty)
\times \partial M\,.$$
{\it Then the \Rm $g'$ cannot be conformally deformed in a complete \Rm of
uniformly positive scalar curvature outside a compact set.}\\[0.1in]
{\bf Proof.} \ \ For the metric $g' = dt^2 + f^2 (t) g$ the Laplacian is 
given
by 
$$
\Delta_{g'} \,u = {{\partial^2 u}\over {\partial t^2}} + {{n f ' 
(t)}\over {f
(t)}} {{\partial u}\over {\partial t}} + {1\over {f^2 (t)}}
\Delta_g u\,, \leqno (1.13)
$$
where $\Delta_g$ is the Laplacian for the \Rm $g$. Let $g'' = u^{4/(n - 
1)} g'$
be a conformal deformation of $g'$, where
$u$ is a positive smooth function. Assume that $g''$ is complete and the
scalar curvature of $g''$, $R_{g''} \ge c^2$ outside a compact set, where $c$
is positive constant. Using (1.2) and the conformal scalar curvature eqaution
[7], we have
$$
{{\partial^2 u}\over {\partial t^2}} + {{n f' (t)}\over {f (t)}} {{\partial
u}\over {\partial t}} + {1\over {f^2 (t)}}
\Delta_g u - {1\over {f^2 (t)}} {{n - 1}\over {4n}}  \{2nf (t) f'' (t) 
+ n (n - 1) |f' (t)|^2 \} u \le - c^2 {{(n - 1)}\over {4n}} u^{{n + 
3}\over {n -
1}}
$$ 
for all $t \ge t'$, where $t' > 2$ is a constant. Fix a value $t
\ge t'$ and integrate the above inequality with respect to the fixed 
Riemannian
manifold $(N, g)$, we have
\begin{eqnarray*} 
(1.14)& \ & {{\partial^2}\over {\partial t^2}} (\int_N u (x, t)
dv_g) + {{n f ' (t)}\over {f (t)}} {{\partial}\over {\partial t}} (\int_N 
u (x,
t) dv_g)\\ & - &{1\over {f^2 (t)}} {{n - 1}\over {4n}}  \{ 2nf (t)
f'' (t)  + n (n - 1) |f' (t)|^2
\}  (\int_N u (x, t) dv_g) \le - c^2 {{(n - 1)}\over {4n}} \int_N u^{{n +
3}\over {n - 1}} dv_g\,,
\end{eqnarray*}  
where we have used Green's identity. For $t \ge t'$, let 
$$U (t) = \int_N u (x, t) dv_g\,.$$
Using the H\"older inequality we have
$$\int_N u^{{n + 3}\over {n - 1}} dv_g \ge {{ [\int_N u dv_g]^{{n + 3}\over
{n - 1}} }\over { [{\mbox{Vol}} (N, g)]^{4\over {n + 3}} }}\,.$$
Thus we have
$${{d^2 U}\over {d t^2}} + {{n f' (t)}\over {f (t)}} {{d
U}\over {d t}} - {1\over {f^2 (t)}} {{n - 1}\over
{4n}}  \{2nf (t) f'' (t)  + n (n - 1) |f' (t)|^2 \} U \le - \epsilon^2 
U^{{n +
3}\over {n - 1}}\,,
\leqno (1.15)$$
where $\epsilon$ is a positive constant given by 
$$\epsilon^2 = {{c^2(n - 1)}\over {4n [{\mbox{Vol}} (N, g)]^{4\over {n +
3}}}}\,.$$ 
In this case $f (t) = t^{1/(n + 1)}$ for $t > 2$, we have 
$$
U'' + {n\over {n + 1}} {{U'}\over t} - {{n - 1}\over {4 (n + 1)}} {{U}\over
{t^2}} \le - \epsilon^2 U^{{n + 3}\over {n - 1}}\,. \leqno (1.16)
$$ 
Let $U (t) = t^\alpha v (t)$ for $t > t'$, where $\alpha$ is a constant 
to be
determined later and $v (t)$ is a positive smooth function. We have
$$U'' (t) = \alpha (\alpha -1) t^{\alpha - 2} v (t) + 2 \alpha t^{\alpha 
- 1}
v' (t) + t^\alpha v'' (t)\,.$$
Equation (1.16) gives
$$t^{\alpha - 2} v [\alpha^2 - {1\over {n + 1}} \alpha - {{n -
1}\over {4 (n + 1)}} ] + t^{\alpha - 1} v' [{n\over {n + 1}} + 2\alpha] + 
t^\alpha v'' (t) \le - \epsilon^2 t^{ {{n + 3}\over {n - 1}} \alpha} 
v^{{n +
3}\over {n - 1}}\,. \leqno (1.17)$$
The quadratic term
$$q (\alpha) = \alpha^2 - {1\over {n + 1}} \alpha - {{n - 1}\over {4 (n + 
1)}}$$
has zeros at
$${1\over 2} ({1\over {n + 1}} \pm \sqrt{ {1\over {(n + 1)^2}} + {{n -
1}\over {n + 1}} } )\,.$$
For 
$$\alpha < {1\over 2} ({1\over {n + 1}} - \sqrt{ {1\over {(n + 1)^2}} + 
{{n -
1}\over {n + 1}} } )\,,$$
we have $q (\alpha) > 0$. If we take $\alpha = - (n - 1)/2$, then
$$\alpha - 2 - {{n + 3}\over {n - 1}} \alpha = 0\,.$$  
From (1.17) we have 
$$t v' ({n \over {n + 1}} + 2\alpha) + t^2 v'' \le - \epsilon^2 v^{{n + 
3}\over
{n - 1}}\,.$$ 
Or
$$(t^{ {n\over {n + 1}} + 2\alpha } v' )' \le - \epsilon^2 t^{ {n\over {n 
+ 1}} +
2\alpha - 2} v^{{n + 3}\over {n - 1}}\,.\leqno (1.18)$$
Integrating both sides of (1.18) we obtain
$$t^{ {n\over {n + 1}} + 2\alpha } v' (t) \le  t_o^{ {n\over {n + 1}} + 
2\alpha
} v' (t_o) - \epsilon^2 \int_{t_o}^t  s^{ {n\over {n + 1}} + 2\alpha - 2} v
(s)^{{n + 3}\over {n - 1}} ds\,.
\leqno (1.19)$$ 
If there exists a $t_o > t'$ such that $v' (t_o) < 0$, then we have
$$
t^{ {n\over {n + 1}} + 2\alpha } v' (t)  \le -c^2
$$ 
for all $t > t_o$, where $c$ is a positive constant. Thus
$$v' (t) \le -c^2 t^{n - 1 - {n\over {n + 1}} }$$
for all $t > t_o\,.$ Therefore we have $v (t) \le 0$ for $t$ large, which 
is a
contracdiction. Hence $v' (t) \ge 0$ for all $t > t'$. Then the inequality
in (1.19) gives
$$
\tau^{ {n\over {n + 1}} + 2\alpha } v' (\tau) \ge \epsilon^2 
\int_{\tau}^t  s^{
{n\over {n + 1}} + 2\alpha - 2} v (s)^{{n + 3}\over {n - 1}} ds\,, \leqno
(1.20)
$$ 
where $t > \tau > t'$. As $v' (t) \ge 0$, we have
$$\tau^{ {n\over {n + 1}} + 2\alpha } v' (\tau) \ge \epsilon^2 v 
(\tau)^{{n +
3}\over {n - 1}} \int_{\tau}^t  s^{ {n\over {n + 1}} + 2\alpha - 2} ds\,.$$
As $\alpha = -(n - 1)/2$ and $n \ge 3$, we have
$${n\over {n + 1}} + 2\alpha - 2 < - 1\,.$$
Therefore after integration and let $t \to \infty$, we have
$$v' (\tau) \ge c_1 ({1\over \tau}) v (\tau)^{{n + 3}\over {n - 1}} \ge 
c_2 {{v
(\tau)}\over \tau}$$ for all $\tau > 2$, where $c_1$ and $c_2$ are positive
constants. We have made use of the fact that $v' (t) > 0$ for all $t > 2$
implies that $v (t)$ is bounded from below by a positive constant. Thus
$$v (t) \ge C t^{c_2}$$
for all $t > 3$, where $C$ is a positive constant. Substitute into the right
hand side of (1.20) gives
$$\tau^{ {n\over {n + 1}} + 2\alpha } v' (\tau) \ge C_1 \int_{\tau}^t  s^{
{n\over {n + 1}} + 2\alpha - 2 + c_2 {{n + 3}\over {n - 1}}}  ds\,,$$ 
for $t > \tau > 3$. As in above, after integration we obtain
$$v' (\tau) \ge c_3 {{v (\tau)}\over {\tau^\delta }} \ge C_2 {{v (\tau)}\over
\tau}$$
where $\delta < 1$ is a positive constant, $c_3$ and $C_2$ are positive
constant. Furthermore, we may assume that 
$$C_1 \ge {n\over {n + 1}} + 2\alpha - 2$$
for all $t$ large enough. Thus 
$$v (t) \ge C_3 t^{C_2}$$ 
for all $t$ large enough, where $C_3$ is a positive constant. Substitute into
the first inequality in (1.20) gives
$$\tau^{ {n\over {n + 1}} + 2\alpha } v' (\tau) \ge C_4 \int_{\tau}^t ds 
\leqno
(1.21)$$ 
for all $\tau$ large enough, where $C_4$ is a positive constant. But
this is impossible as the right hand side of (1.21) tends to infinity as $t
\to \infty$.\qed 

\vspace{0.3in}

{\bf \Large 3. \ \ Boundary components in class (C)}

\vspace{0.3in}

In this section we assume that at least one of the boundary components of
$\overline M$, namely $N_1$, belongs to class (C). Then any \Rm $g$ on $N_1$
would have the scalar curvature negative somewhere. Take a \Rm $g_1$ on
$N_1$ with $R (g_1) = -n (n - 1)$. Then equations (1.3) becomes
$${{4n}\over {n + 1}} u'' + n (n - 1) u^{{n - 3}\over {n + 1}} + R u = 0\,.
\leqno (3.1)$$
{\bf Lemma 3.2.} \ \ {\it Assume that $R \in C^\infty ([2, \infty))$ is a
negative function such that $R \ge - a^2$ for some positive constant $a$ and}
$$R (t) \le - {C\over {t^\alpha}} \ \ \ \ {\mbox{for}} \ \ \ t \ge t_o\,,$$
{\it where $t_o > 2$, $C$ and $\alpha \le 2$ are positive constants. If 
$\alpha =
2$, we assume that $C > n (n - 1)$. Then equation (3.1) has a positive
solution on $(2, \infty)$.}\\[0.1in]
{\bf Proof.} \ \ If $\alpha < 2$, then we let $u_+ = c_+ + t^m$ and $u_- 
= c_-$,
where $c_+, c_-$ and $m$ are positive numbers. If we take $c_+$ and $m$ large
enough and take $c_-$ small, then we have 
$${{4n}\over {n + 1}} u_{+}" + n (n - 1) u_{+}^{{n - 3}\over {n + 1}} + R 
u_{+} 
\le  0\,,$$
$${{4n}\over {n + 1}} u_{-}" + n (n - 1) u_{-}^{{n - 3}\over {n + 1}} + R 
u_{-} 
\ge  0\,.$$ By the upper and lower
solution method, we obtain a positive solution (c.f. [7]). In case 
$\alpha = 2$
and
$C > n (n - 1)$, we may take $u_+ = C_+ t^{(n + 1)/2}$, where $C_+$ is a 
positive
constant. Then
\begin{eqnarray*}
& \ &{{4n}\over {n + 1}} u_{+}'' + n (n - 1) u_{+}^{{n - 3}\over {n + 1}} 
+ R
u_{+}\\
&\le & C_+n (n - 1) t^{{n -3}\over 2} + n (n + 1) C_+^{{n - 3}\over {n 
+1}} t^{{n
-3}\over 2} - C_+ [n (n - 1) + \epsilon] t^{{n -3}\over 2} \le 0\,,
\end{eqnarray*} 
if we take
$C_+$ to be large enough. Here $\epsilon = C - n (n - 1) > 0$ is a positive
constant. Take $u_{-}$ to be a small positive constant. In this case, we 
obtain a
positive solution as in above.\qed 
\hspace*{0.5in}In the above lemma, when $\alpha = 2$ and $C \ge n
(n - 1)$, we have the following.\\[0.1in]  
{\bf Lemma 3.3.} \ \
{\it Suppose that $N_1$ belongs to class (C). Let $g$ be a \Rm on $N_1$. 
On the
end
$(2,
\infty)
\times N_1$, there does not exist a warped product metric}
$$g' = dt^2 + f^2 (t) g$$
{\it with scalar curvature}
$$R \ge -{{n (n - 1)}\over {t^2}}$$
{\it for all $x \in N_1$ and $t > t_o > 2$, where $t_o$ is a 
constant.}\\[0.1in]
{\bf Proof.} \ \ Assume that we can find a warped product metric on $(2, 
\infty)
\times N_1$ with 
$$R \ge -{{n (n - 1)}\over {t^2}}$$
for all $x \in N_1$ and $t > t_o > 2$. We may assume that the scalar 
curvature of
$g$ is equal to
$-\kappa^2$ at $x_o \in N_1$, where $\kappa$ is a positive constant. With 
$u (t)
= f^{{n + 1}\over 2} (t)$ and at
$x_o$, by (1.3) we have
$${{4n}\over {n + 1}} {{u''}\over {u}} + {{\kappa^2}\over {u^{4\over {n + 
1}} }}
\le {{n (n - 1)} \over {t^2}}\,. \leqno (3.4)$$ In particular
$${{u''}\over {u}} \le  {{(n +
1) (n - 1)} \over {4 t^2}}\,.$$ That is,
$$t^{{n + 1}\over 2} u'' \le {{(n + 1) (n - 1)} \over {4}} t^{{n - 3}\over
2}u\,.$$  Consider the inequality
$${{u'' (t)}\over {u (t)}} \le {C \over {t^2}}$$ for $t > t_o > 2$, where 
$C \ge
1$ is a constant. Let $\varepsilon > 1$ be a constant such that $\varepsilon
(\varepsilon - 1) = C$. Then we have
$$
t^\varepsilon u'' (t) \le \varepsilon (\varepsilon - 1) t^{\varepsilon - 
2} u
(t)\,.
$$  Upon integration from
$t_1 \ge t_o$ to $t > t_1$, and using integration by parts twice, we obtain
$$t^\varepsilon u' (t) - \varepsilon t^{\varepsilon - 1} u (t)  -
t_1^\varepsilon u' (t_1) + \varepsilon t_1^{\varepsilon - 1} u (t_1) +
\varepsilon (\varepsilon - 1) \int_{t_1}^t s^{\varepsilon - 2} u (s) ds 
\le C
\int_{t_1}^t s^{\varepsilon - 2} u (s) ds\,.$$   Therefore we have
$$t^\varepsilon u' (t) - \varepsilon t^{\varepsilon - 1} u (t) \le 
t_1^\varepsilon u' (t_1) - \varepsilon t_1^{\varepsilon - 1} u (t_1)\,. 
\leqno
(3.5)$$ If there is a number $t_1 \ge t_o$ such that
$u' (t_1)
\le 0$, then we have 
$$t^\varepsilon u' (t) - \varepsilon t^{\varepsilon - 1} u (t) \le 0\,.$$ 
This
gives
$$(\ln u (t) )' \le \varepsilon (\ln t)'\,.$$  Hence
$$u (t) \le c t^\varepsilon$$ for all $t > t_1$, where $c$ is a positive
constant. If $u' (t) > 0$ for all $t
\ge t_o$, then $u (t) \ge c'$ for some positive constant $c'$. Let $C$ be a
positive constant such that 
$$t_1^\varepsilon u' (t_1) -  \varepsilon t_1^{\varepsilon -1} u (t_1)
\le C\,,$$ 
then (3.5) gives 
$$t^\varepsilon u' (t) - \varepsilon t^{\varepsilon - 1} u (t) \le C$$ 
for all
$t > t_1$.  Thus 
$${{u' (t)}\over {u (t)}} \le {\varepsilon \over t} + {C \over {u (t)
t^\varepsilon}}
\le  {\varepsilon \over t} + {C \over {c t^\varepsilon}}\,.$$ Integrating 
from
$t_1$ to $t$ we have
$$\ln {{u (t)}\over {u (t_1)}} \le \varepsilon \ln ({{t}\over {t_1}}) + 
{C \over
{c' t_1^{\varepsilon - 1}}} \le \varepsilon \ln ({{C' t}\over 
{t_1}})\,,$$  as
$\varepsilon > 1$. Here $C'$ is a positive constant such that $\ln C' \ge 
C/ (c'
t_1)$. Hence we again obtain the inequality
$$u (t) \le b t^\varepsilon$$ for some positive constant $b$ and for all 
$t \ge
t_1$. Thus we find a constant $c > 0$ such that 
$$u (t) \le c t^\varepsilon \leqno (3.6)$$ 
for all $t \ge t_1$. In case $C = (n +
1) (n - 1)/4$ we take $\varepsilon = (n + 1)/2$, we have
$$u (t) \le c t^{{n + 1}\over 2}\,.$$ Then 
$$
{{\kappa^2}\over {u^{4 \over {n + 1}} }} \ge {c' \over {t^2}}\,,
$$ 
where $c'$ is a positive constant. Hence
(3.4) gives
$${{u''}\over {u}} \le  {{(n + 1) (n - 1) - \delta} \over {4 t^2}}\,,$$ where
$\delta > 0$ is a costant. Similarly we have
$$u (t) \le c t^{ {{n + 1}\over 2} - \delta '}$$ and
$${{\kappa^2}\over {u^{4 \over {n + 1}} }} \ge {c'' \over {t^{2 -
\epsilon}}}\,.$$  for some positive constants $\delta '$, $\epsilon$ and
$c''$. Thus (3.4) gives
$$u'' (t) \le 0$$ for $t$ large and hence $u (t) \le C t$ for some 
contant $C >
0$. From (3.4) we have 
$${{u'' (t) }\over {u (t)}} \le -{{\kappa^2}\over { (Ct)^{4\over {n + 1}} 
}} +
{{n(n - 1)}\over {t^2}} \le  - {c\over t}$$
for $t$ large enough, as $n \ge 3$. Here $c$ is a positive constant. We 
have 
$$u' (t) - u' (t') \le - c \int_{t'}^t {{u (s)}\over s} ds\,,  \ \ \ \ t >
t_1\,.$$  
If  $u' (t') \le
0$ for some $t'$, then $u' (t) \le - c_1$ for some positive constant $c_1$.
Hence $u (t) \le 0$ for $t$ large enough, contradicting the fact that $u$ is
positive. If
$u' (t) > 0$ for all $t$ large, then 
$$\int_{t'}^t {{u (t)}\over t} dt \ge u (t') \int_{t'}^t {{u (s)}\over s} 
ds\to
\infty\,.$$ Thus $u'$ has to be negative for some $t$ large.\qed 
\hspace*{0.5in}The result in lemma 3.3 is
almost sharp as we can get as close to $-n (n - 1)/t^2$ as possible. For
example, Let $R (g) = - n (n - 1)$ let $f (t) = t \ln t$ for $t> 2$. We have
$$
R = -{1\over {t^2}} [ n (n - 1) {{(\ln t + 1)^2}\over {(\ln t)^2}} +
{{2n}\over {\ln t}} + {{n (n - 1)}\over {(\ln t)^2}} ]\,.\leqno (3.7)
$$  
We can show that it
is not possible to conformally deform the metric 
$$
dt^2 + (t \ln t)^2 g
$$ 
into a complete metric of nonnegative scalar
curvature outside a compact set. Acturally we can show more.  We first 
note that
if 
$$f (t) f'' (t) \le -c^2$$ for some positive constant $c$ and for $t > 
t_o > 2$,
then
$$f' (t) \le f (t') - \int_{t'}^t {{c^2}\over {f (s)}} ds\,,$$ where $t > 
t' >
t_o$. As $f'' (t) \le 0$, we have $f (t) \le C t$ for all $t > t_o$, 
where $C$
is a positive constant. Therefore 
$$
f' (t) \le f (t') - ({c\over C})^2 \int_{t'}^t {1\over t} ds \to -\infty 
\ \
\ \ {\mbox{as}} \ \ t \to \infty\,.
$$ 
Thus $f (t) \le 0$ when $t$ is large. As in theorem 1.12, we assume
that the boundary of $\overline M$ is connected. And the case where $\partial
\overline M$ has more than one connected component is similar.\\[0.1in]

{\bf Theorem 3.8.} \ \ {\it For $n \ge 3$, let $\overline M$ be a compact
$(n + 1)$-manifold with boundary $N$. Suppose that $N$ belongs to class
(C) and $g$ is a \Rm on $N$ with scalar curvature $R (g) \le - \kappa^2$ for
some positive constant $\kappa$. Let $g'$ be a complete \Rm on $M$ with} 
$$
g' (t, x) = dt^2 + f (t) g (x) \ \ \ \ on \ \ \ (2,
\infty)
\times N\,.
$$
{\it Suppose that $f (t) f'' (t) \ge -c^2$ for all $t$ large, where $c^2 =
2(\kappa^2 -\delta)/(3n + 1)$ and $\delta$ is constant such that $0 <
\delta < \kappa^2$. Assume that either (i)
$f (t) < C t
\ln t$ or (ii) $f (t) \ge C t^\alpha$ with $\alpha > 1$, for $t$ large, then
there does not exist a complete \Rm conformal to
$g'$ with nonnegative scalar curvature outside a compact set.}\\[0.1in]  
{\bf Proof.} \ \ Let $g'' = u^{4/(n - 1)} g'$ is a \Rm conformal to $g'$, 
where
$u$ is a positive smooth function. Assume that $g''$ is complete and
the scalar curvature of $g''$ satisfies $R_{g''}
\ge 0$ outside a compact set. The conformal scalar curvature equation [7] 
gives 
$${{\partial^2 u}\over {\partial t^2}} + {{n f' (t)}\over {f (t)}} {{\partial
u}\over {\partial t}} + {1\over {f^2 (t)}}
\Delta_g u - {1\over {f^2 (t)}} {{n - 1}\over {4n}}  \{ R (g) - 2nf (t) 
f'' (t) 
- n (n - 1) |f' (t)|^2 \} u \le 0\,.$$ for all $t  \ge t'$. Fix a value 
$t \ge
t'$ and integrate the above inequality with respect to the fixed Riemannian
manifold $(N, g)$, we have
\begin{eqnarray*} & \ & {{\partial^2}\over {\partial^2 t}} (\int_N u (x, t)
dv_g) + {{n f ' (t)}\over {f (t)}} {{\partial}\over {\partial t}} (\int_N 
u (x,
t) dv_g)\\ & + &{1\over {f^2 (t)}} {{n - 1}\over {4n}}  \{ \kappa^2 + 2nf (t)
f'' (t)  + n (n - 1) |f' (t)|^2
\}  (\int_N u (x, t) dv_g) \le 0\,,
\end{eqnarray*}  where we have used Green's theorem and the fact that $R_g
\le -\kappa^2$ for some positive constant $\kappa$. For $t \ge t'$, let 
$$U (t) = \int_N u (x, t) dv_g\,.$$ Then we have 
$$
U'' +  {{n f ' (t)}\over {f (t)}} U' + {1\over {f^2 (t)}} {{n - 1}\over 
{4n}} 
\{ \kappa^2 + 2nf (t) f'' (t)  + n (n - 1) |f' (t)|^2
\} U \le 0\,.\leqno (3.9)
$$  
For $t \ge t' > 0$, let $U (t) = f^\alpha (t) v (t)$, where
$\alpha$ is a constant to be determined later. Then we have
\begin{eqnarray*} U' (t) & = &\alpha f^{\alpha - 1} (t) f' (t) v (t) + 
f^\alpha
(t) v' (t)\,,\\  U'' (t) & = & \alpha (\alpha - 1) f^{\alpha - 2} (t) |f' 
(t)|^2
v (t) + \alpha f^{\alpha - 1} (t) f'' (t) v (t) + 2
\alpha f ^{\alpha - 1} (t) f' (t) v' (t) + f^\alpha (t) v'' (t)\,.
\end{eqnarray*} Substitute into (3.9) we have
\begin{eqnarray*} (3.10) \ \ \ & \ & f^{\alpha - 2} (t) |f' (t)|^2 v (t) 
[\alpha
(\alpha - 1) + n \alpha + {{(n - 1)^2}\over 4}] + f^\alpha (t) v'' (t)\\ 
& \ &
\ \ \  + (n + 2\alpha) f^{\alpha - 1} (t) f' (t) v' (t) +  f^{\alpha - 2} 
(t) v
(t) [\kappa^2 + (2n + \alpha) f (t) f'' (t) ]\le 0
\end{eqnarray*} for all $t \ge t' > 0$. As 
$$
\alpha (\alpha - 1) + n \alpha + {{(n - 1)^2}\over 4} = [\alpha + {{n - 
1}\over
2}]^2 \ge 0
$$  
and $f (t)f'' (t) \ge 2(-\kappa^2 + \delta)/(3n + 1)$, we have
$$f^\alpha (t) v'' (t) + (n + 2\alpha) f^{\alpha - 1} (t) f' (t) v' (t) \le
-[\kappa^2 + (2n + \alpha) {{2(-\kappa^2 + \delta)}\over {(3n + 1)}}] 
f^{\alpha -
2} (t) v (t)\,.$$ Or
$$
(f^{n + 2\alpha} (t) v' (t))' \le -[\kappa^2 + (2n + \alpha) 
{{2(-\kappa^2 +
\delta)}\over {(3n + 1)}}] f^{n + 2 \alpha - 2} (t) v (t)\,.
\leqno (3.11)
$$
Choose $\alpha$ such that $n + 2 \alpha = 1$, that is, $\alpha = - (n - 
1)/2$, we
have
$$(f (t) v' (t) )' \le -\delta {{v (t)}\over {f (t)}}\,.$$ Upon 
integration we
have
$$
f (t) v' (t) - f (t_o) v' (t_o) \le -\delta  \int_{t_o}^t {{v (s)}\over {f
(s)}} ds\,,\leqno (3.12)
$$ where $t > t_o \ge t'$. If $v' (t) > 0$ for all $t \ge t'$, then
we have
$$
f (t) v' (t) \le f (t_o) v' (t_o) - \delta v (t_o) \int_{t_o}^t {1\over {f
(s)}} ds\,.
$$ 
As $f (t) \le C t \ln t$, for some positive constants $C$ and for $t > t_o$,
we have
$$\int_{t_o}^t {1\over {f (s)}} ds \ge {1\over C} \int_{t_o}^t {1\over {s 
\ln s}}
ds = {1\over C} (\ln
\ln t - \ln \ln t_o) \to \infty$$ as $t \to \infty$. That is, $v' (t) < 
0$ when
$t$ is large. Thus we can find a
$t_o$ such that $v' (t_o) \le 0$. Hence
$$f (t) v' (t) \le - c'^2 \int_{t_o}^t {{v (s)}\over {s \ln s}} ds\,,$$ 
that is
$v' (t) \le 0$ for all $t \ge t_o$. Here $c'$ is a positive constant. We have
$$
f (t)  v' (t) \le - c'^2 v (t) \int_{t_o}^t {1\over {s \ln s}} ds\,. \leqno
(3.13) 
$$ For any
positive constant $C > 0$, we can find $t'' > t_o$ such that for all $t
\ge t''$, (3.13) gives 
$$(\ln v (t) )' \le - c'^2 {{\ln \ln t}\over {t \ln t}}\,.$$ Therefore 
$$\ln {{v (t)}\over {v (t_o)}} \le -c''^2 (\ln \ln t)^2\,,$$ 
where $c''$ is a positive constant. Thus 
$$v (t) \le {{C'}\over {(\ln t)^\beta}}$$ for some positive constants 
$C'$ and
$\beta$ and for $t > t_o$.  Thus 
$$
U (t) \le {{C'}\over {f^{{n - 1}\over 2} (t) (\ln t)^\beta}} \le {{C''} \over
{(t \ln t)^{{n - 1}\over 2}(\ln t)^\beta}}\,.\leqno (3.14)
$$ 
Since $n \ge 3$, we have
$2/(n - 1) \le 1$. By H\"older's inequality, we have
$$\int_N u^{2 \over {n - 1}} dv_g \le [ {\mbox {Vol}} (N, g) ]^{{n - 
3}\over {n
- 1}} (\int_N u dv_g)^{2\over {n - 1}} \le C'' U^{2\over {n - 1}} (t) \le
{{C'}\over {t (\ln t)^\gamma}}\,,$$ where $\gamma > 1$ is a positive 
constant.
Thus
$$
\int_{t''}^\infty \int_N u^{2\over {n - 1}} dv_g dt < C' \int_{t''}^\infty
{{dt}\over {t (\ln t)^\gamma}} < \infty\,.
\leqno (3.15)
$$ 
Hence we can
find $x_o \in N$ such that 
$$\int_{t''}^\infty u^{2\over {n - 1}} (x_o, t) dt < \infty\,.$$ Therefore
the curve $\gamma (t) = (t, x_o)$ for $t \in (t''\,, \infty)$ has finite
length in the \Rm $g'' = u^{4/(n - 1)} g'\,,$ that is, the metric
$g''$ is not complete. In case
$f (t) > C t^\alpha$ for
$\alpha > 1$, then (3.12) gives
$$f (t) v' (t) \le C'$$ or
$$v' (t) \le {{C'}\over {t^\alpha}}$$ Thus $v (t) \le C''$ for all $t > 
t_o$ and
hence 
$$U (t) \le {C\over {t^{ {{(n - 1)}\over 2}\alpha} }}\,.$$ As $\alpha > 
1$, we
can conclude as above that the metric $g''$ is not complete.\qed
\hspace*{0.5in}We note that functions of the type $f (t) = C t^\alpha (\ln
t)^\beta$ satisfy the conditon of theorem 3.8, where 
$\alpha, \beta$, and $C > 0$ are constants, where if $\alpha = 1$, then
$\beta \le 1$.

\vspace{0.3in}

{\bf \Large 4. \ \ Scalar curvature of polar type \Rms}

\vspace{0.3in}

Given a compact $n$-manifold $N$ and a constant $a > 0$, consider the 
following
metric on $(a,
\infty)
\times N$:
$$\og (t, x) = dt^2 + f^2 (t, x) g (x)\,, \leqno (4.1)$$
where
$f$ is a positive smooth function on $(a, \infty) \times N$ and $g$ is a
fixed \Rm on $N$. Let $R (f^2 (t, \bullet) g )$ be
the scalar curvature on $N$ corresponding to the \Rm $f^2 (t, \bullet) 
g$, that
is, the \Rm conformal to $g$ with conformal factor $f^2 (t, \bullet)$, 
where $t$
is treated as a constant. Let $\overline R$ be the scalar curvature of 
the \Rm
$\og$. We show that (c.f. [3])
$$
\overline R (t, x) = R (f^2 (t, \bullet) g ) (t,
x) - {1\over {f^2 (t, x)}} [ 2n f (t, x) {{\partial^2 f}\over {\partial t^2}}
(t, x) + n (n - 1) |{{\partial f}\over {\partial t}} (t, x)|^2 ]\,. \leqno
(4.2)
$$ 
In particular, if $f$ depends on $t$ only, then 
$$
R (f^2 (t, \bullet) g ) (t, x) = {1\over {f^2 (t)}} R (g)\,,
$$
which gives equation (1.2). In general, if $n \ge 3$ and we let $\mu
(\bullet) = f^{{n - 2}\over 2} (t, \bullet)$, then
$$ R (f^2 (t, \bullet) g ) = c_n^{-1} \mu^{- {{n + 2}\over {n - 2}}
}[c_n R (g) \mu - \Delta_g
\mu]\,, \leqno (4.3)$$
where $c_n = (n - 2)/[4 (n - 1)]\,$ [7].\bk
Let $(x_1, x_3,..., x_n)$ be a set of local
coordinates for $N$ and let $x_o = t$. Since $R_{0jk0} = R_{j00k}\,,$ we have
\begin{eqnarray*}
(4.4) \ \ \ \ \ \ \ \ \overline R & = & \sum_{0 \le i, j, k, l \le n}
\og^{jl}\og^{ik}
\oR_{ijkl}\\ & = & 2 \sum_{1 \le j, k \le n} \og^{jl}
\oR_{0j0l} +
\sum_{1
\le i, j, k, l \le n}
\og^{jl}\og^{ik} \oR_{ijkl}\,, \ \ \ \ \ \ \ \ \ 
\end{eqnarray*}
where $\oR_{ijkl}$ is the Riemannian curvature tensor for the metric 
$\og$ and
is given by
$$\oR_{ijkl} = {1\over 2} 
(   {{\partial^2 \og_{jk}}\over
{\partial x_i \partial x_l}} 
+ {{\partial^2 \og_{il}}\over {\partial x_j \partial
x_k}} - {{\partial^2 \og_{ik}}\over {\partial x_j \partial x_l}} 
- {{\partial^2
\og_{jl}}\over {\partial x_i \partial x_k}}) 
+ \sum_{0 \le a,
b \le n} \og_{ab} (\Gamma^b_{il}\Gamma^a_{jk} 
-\Gamma^b_{ik}\Gamma^a_{jl})\,.$$
Here 
$$\Gamma^i_{jk} = {1 \over 2} \sum_{a= 0}^n \og^{ia} ({{\partial \og_{aj} 
}\over
{\partial x_k}} + {{\partial \og_{ak} }\over {\partial x_j}} - {{\partial
\og_{ik} }\over {\partial x_a}} )\,, \ \ \ 0 \le i, j, k \le n$$
are the Christoffel symbols for $\og$. Using 
\begin{eqnarray*}
(4.5) \ \ \ \ \ \ \  \ \Gamma^0_{jk} & = & - f {{\partial f}\over 
{\partial t}}
g_{jk}\,, \ \ \ \ 1 \le i, j \le n;\\
\Gamma^0_{0k} & = & 0\,, \ \ \ \ 0 \le k \le n;\\ 
\Gamma^i_{00} & = & 0\,, \ \ \ \ 0 \le i \le n;\\
\Gamma^i_{0k} & = & {1\over f} {{\partial f}\over 
{\partial t}} \delta_{ik}\,, \ \ \ \ 1 \le i,  k \le n;\\
\Gamma^i_{jk} & = &  {1 \over 2} \sum_{a= 1}^n \og^{ia} ({{\partial \og_{aj}
}\over {\partial x_k}} + {{\partial \og_{ak} }\over {\partial x_j}} - 
{{\partial
\og_{ik} }\over {\partial x_a}} )\,, \ \ \ 1 \le i, j, k \le n\,, 
\end{eqnarray*}
we obtain $\oR_{0000} = 0$ and for $1 \le j, k \le n$, we have
$$\sum_{1 \le j, k \le n} \og^{jl} \oR_{0j0l} = -{n \over {f}} {{\partial^2
f}\over {\partial t^2}}\,. \leqno (4.6)$$ 
For $1 \le i,j,k,l \le n$, we have
\begin{eqnarray*}
\oR_{ijkl} & = & {1\over 2}  (  {{\partial^2 \og_{jk}}\over
{\partial x_i \partial x_l}} + {{\partial^2 \og_{il}}\over {\partial x_j 
\partial
x_k}} - {{\partial^2 \og_{ik}}\over {\partial x_j \partial x_l}} - 
{{\partial^2
\og_{jl}}\over {\partial x_i \partial x_k}}) + \sum_{1 \le a, b \le n} 
\og_{ab}
(\Gamma^b_{il}\Gamma^a_{jk} -
\Gamma^b_{ik}\Gamma^a_{jl})\\
& \ \ & \ \ + \og_{00} (\Gamma^0_{il}\Gamma^0_{jk} -
\Gamma^0_{ik}\Gamma^0_{jl})
\end{eqnarray*}
The first two terms on the right hand side of the above equation has no
derivatives with respect to $t$ and by (4.5) it is equal to $R (f^2 (t, 
\bullet)
g )_{ijkl}\,,$ the Riemannian curvature tensor for the \Rm $f^2 (t, 
\bullet) g$, 
for $1 \le i,j,k,l\,.$ And 
$$(\Gamma^0_{il}\Gamma^0_{jk} -
\Gamma^0_{ik}\Gamma^0_{jl}) = f^2 |{{\partial f}\over {\partial t}}|^2
(g_{il}g_{jk} - g_{ik}g_{jl})\,.$$
Thus we have
\begin{eqnarray*}
(4.7) \ \ \ \ \sum_{1 \le i, j, k, l \le n}
\og^{jl}\og^{ik} \oR_{ijkl} & = & \sum_{1 \le i, j, k, l \le n} 
\og^{jl}\og^{ik}
R (f^2 (t, \bullet) g )_{ijkl}\\
& \ & \ \ \  + \sum_{1 \le i, j, k, l \le n} \og^{jl}\og^{ik}
f^2 |{{\partial f}\over {\partial t}}|^2 (g_{il}g_{jk} - g_{ik}g_{jl})\\ 
& = & R (f^2 (t, \bullet) g ) - {{n (n - 1)}\over {f^2}}|{{\partial f}\over
{\partial t}}|^2\,.
\end{eqnarray*}
Substitute (4.6) and (4.7) into (4.4) we have obtained the desired
formular. Using (4.4), the curve $\gamma (t) = (t, x)$ is a geodesic for 
$t >
2$,  where $x \in N$. The metric $\og$ as defined in (4.1) is of polar 
type (c.f.
[3]).   
\\[0.1in]  
{\bf Theorem 4.8.} \ \ {\it For $n \ge 3$, let $g$ be a \Rm on $N$ with
nonpositive total scalar curvature. Then the
scalar curvature $\oR$ of the \Rm $\og$ in (4.1) cannot be uniformly 
positive for
$t$ large enough.}\\[0.1in] 
{\bf Proof.} \ \ Assume that $\oR \ge \alpha^2$ for all
$t > t'$ and $x \in N$, where $\alpha$ and $t' > a$ are positive 
constants. From
(4.3) we have
$$ 
R (f^2 (t, \bullet) g ) = {1\over f^2} [ R (g)
-  c_n^{-1} {{ \Delta_g \mu}\over \mu}]\,.
$$ 
Substitute into equation (4.2) we have   
\begin{eqnarray*}
(4.9) \ \ \ \ \ \ \ f^2 \overline R (t, x) & = & -  c_n^{-1} {{ \Delta_g
\mu}\over \mu}  + R (g) -
2n f (t, x) {{\partial^2 f}\over {\partial t^2}} (t, x) - n (n - 1) 
|{{\partial
f}\over {\partial t}} (t, x)|^2 \\
& = & - c_n^{-1} {{ \Delta_g \mu}\over \mu} - n {{\partial^2 (f^2)}\over
{\partial t^2}} - n (n - 3) |{{\partial f}\over {\partial t}}|^2\,,
\end{eqnarray*}
where we have used the formular
$$
-2n f {{\partial^2 f}\over {\partial t^2}} = -n {{\partial^2 (f^2)}\over
{\partial t^2}} + 2n |{{\partial f}\over {\partial t}}|^2\,.
$$
Fix a $t > t' $ and integrate both sides of (4.9) with
respect to $(N, g)$, we have
$$
\int_N f^2 \oR dv_g = - c_n^{-1} \int_N {{|\bigtriangledown \mu|^2}\over
{\mu^2}} dv_g + \int_N R (g) dv_g - n 
\int_N {{\partial^2 (f^2)}\over {\partial t^2}} dv_g - n (n - 3) \int_N
|{{\partial f}\over {\partial t}}|^2 dv_g\,,
\leqno (4.10)
$$
where we have used Green's identity. Using the fact that $n \ge 3$, $\oR
\ge \alpha^2$ for $t > t'$ and 
$$
\int_N R (g) dv_g \le 0\,,
$$
we have
$$
a^2 F \le -n {{d^2 F}\over {d t^2}}\,,
$$ 
where
$$
F (t) = \int_N f^2 dv_g\,.
$$
Or
$$
F'' (t) \le - b^2 F (t)  \ \ \ \ {\mbox{for}} \ \ \ t > t'\,,
\leqno (4.11)
$$
where $b = \alpha/\sqrt n$ is a positive constant. If there is $t_o > t'$
such that $F' (t_o) < 0$, then integrating both sides of (4.11) gives
$$
F' (t) \le F' (t_o) - b^2 \int_{t_o}^t F (s) ds \ \ \ \ 
{\mbox{for}} \ \ \ t > t_o\,.
\leqno (4.12)
$$
Thus $F (t) \le 0$ for some
$t$ large enough, contradicting that $f$ is a positive function. Therefore
$F' (t) \ge 0$ for all $t > t'$. Then (4.11) gives
$$
F' (t) \le F' (t_o) - b^2 F (t_o) \int_{t_o}^t ds \to - \infty  \ \ \ \
{\mbox{as}} \ \ \ t \to \infty\,,
$$
which contradicts $F' (t) \ge 0$ for all $t > t'$.\qed
\hspace*{0.5in}Conbining the proof of lemma 1.8 and the above theorem 
4.8, it
can be shown that the scalar curvature $\oR$ of $\og$ cannot decay to 
zero too
slowly. More precisely, $\oR$ cannot be bigger than or equal to 
$${{cn}\over {4t^2}}$$
for all $t > t' > a$ and $x \in N$, where $c > 1$ and $t'$ are constants. For
if 
$$\oR \ge {{cn}\over {4t^2}}\,,$$
then (4.10) gives
$$t^2 F'' + {c\over 4} F \le 0$$
for all $t > t'$. The proof of lemma 1.8 shows that $F (t) = 0$ for some
$t$ large. On the other hand it is known that $F (t) > 0$ for all $t >
a$.\\[0.1in] 
{\bf Theorem 4.13.} \ \ {\it For $n \ge 3$, let $g$ be a \Rm on $N$
with negative total scalar curvature. Then the scalar curvature $\oR$ of 
the \Rm
$\og$ in (4.1) cannot be bigger than $-c/t^2$ for some constant $c <
2n$ and for all $t > t' > 2$ and $x \in N$. In particular, $\oR$
cannot be nonnegative for all $t$ large.}\\[0.1in]   {\bf Proof.} \ \ Assume
that 
$$
\int_N R (g) dv_g \le - b^2
$$
and 
$$\oR (t, x) \ge - {c \over {t^2}} \ \ \ \ {\mbox{for}} \ \ t > t' \ \
{\mbox{and}} \ \ x \in N\,,$$
where $c < 2n$ and $b$ is a positive constant. Then (4.10) gives
$$
-{c \over {t^2}} F (t) \le  -b^2  - n F'' (t)  \ \ \ \ {\mbox{for}} \ \ t >
t'\,,
$$
or 
$$
F'' (t) \le - {{b^2}\over n} + {{c'}\over {t^2}} F (t)  \ \ \ \ 
{\mbox{for}} \
\ t > t'\,,\leqno (4.14)
$$
where $c' = c/n < 2$ is a constant. If $c' \le 0$, then (4.14) shows
that $F (t) \le 0$ for $t$ large enough. So we may assume that $c' > 0$. From
(4.14) we have 
$$
{{F'' (t) }\over {F (t)}} \le {{c'}\over {t^2}} \ \ \ \ {\mbox{for}} \
\ t > t'\,.
$$
It follows from the proof of lemma 3.3 that 
$$
F (t) \le C t^\epsilon \ \ \ \ {\mbox{for}} \
\ t > t_o > t'\,, \leqno (4.15)
$$
where $C$ is a positive constant and $\epsilon > 1$ is a positive constant
such that $\epsilon (\epsilon - 1) = c'$. In particular, $\epsilon < 2$. Then
(4.14) and (4.15) imply that $F'' (t) \le - b^2/(2n)$ for all $t$ large
enough. Thus $F (t) \le 0$ when $t$ is large enough.\qed 
\hspace*{0.5in}The Laplacian for the metric $\,\og (t, x) = dt^2 + f^2 
(t, x) g
(x)\,$ is given by
$$
\Delta_\og u = {{\partial^2 u}\over {\partial t^2}} + {n \over f}{{\partial
f}\over {\partial t}} {{\partial u}\over {\partial t}} + {{n - 2}\over {f^3}}
<\bigtriangledown_g f\,, \bigtriangledown_g u >_g + {1\over {f^2}} \Delta_g
u
\leqno (4.16)
$$
for $u \in C^\infty ( (a\,, \infty) \times N)\,.$ Suppose that $u (t, x)$ is
a positive smooth function. The scalar curvature
$R_c$ of the
\Rm 
$$
u^{4\over {n - 1}} (t, x) \og = u^{4\over {n - 1}} (t, x) [dt^2 + f^2 (t, 
x) g
(x)]
$$
is given by the following:
$$
\Delta_\og \,u - c_{n + 1} R_\og u + c_{n + 1} R_c u^{{n + 3}\over {n - 
1}} =
0\,, \leqno (4.17)
$$
where $c_{n + 1} = (n - 1)/(4n)$.\\[0.1in]
{\bf Theorem 4.18.} \ \ {\it For $n \ge 3$, let $\overline M$ be
a compact $(n + 1)$-manifold with boundary $N$. Suppose that $N$ is connected
and $g$ is a \Rm on $N$. Suppose
that $f$ is a positive smooth function on $(2, \infty) \times N$ such that}
$$
|{{\partial f}\over {\partial t}} | \le C_1 {f\over t} \ \ \ \ 
{\mbox{and}} \ \
\ |{{\partial^2 f}\over {\partial t^2}} | \le C_2 {f\over {t^2}}\,,\leqno 
(4.19)
$$
{\it for all $t > t' > 0$ and $x \in N$. Let $\og$ be a complete \Rm on 
$M$ such
that $\,\og (t, x) = dt^2 + f^2 (t, x) g (x)\,$ on $(2, \infty) \times N$.
Assume that $\og$ has scalar curvature
$R_\og
\le -b^2$ for all $t > t'$ and $x \in N$, where $t' > 2$, $C_1$, $C_2$ 
and $b$
are positive constants. Then for any positive smooth function $u$ on $(a, 
\infty)
\times N$ with}
$$|{{\partial u}\over {\partial t}} | \le C u \ \ \ \ {\mbox{for \ all}} 
\ \ \
t > t_1 \ \ \ {\mbox{and}} \ \ x \in N\,, \leqno (4.20)$$
{\it where $t_1 > 2$ and $C$ are positive constant, if $u^{4\over {n - 1}}
\,\og$ is a complete \Rm, then the scalar curvature
$R_c$ of the conformal metric $u^{4\over {n - 1}} \,\og$ cannot not be
nonnegative for all $t$ large\,.}\\[0.1in]
{\bf Proof.} \ \ Assume that the scalar curvature of the \Rm $u^{4\over 
{n - 1}}
\,\og$ is nonnegative for all $t > t_2$ and $x \in N$, where $t_2 > 2$ is a
constant. If $t_o = \max \{t'\,, t_1\,, t_2\,,\}\,,$ then for $t > t_o$,
the scalar curvature equation (4.17) togather with (4.16) and (4.17) give
$$
{{\partial^2 u}\over {\partial t^2}} + {n \over f}{{\partial f}\over 
{\partial
t}} {{\partial u}\over {\partial t}} + {{n - 2}\over {f^3}} 
<\bigtriangledown_g
f\,, \bigtriangledown_g u >_g + {1\over {f^2}} \Delta_g u + c^2 u \le
0\,,\leqno (4.21)
$$
where $c^2 = c_{n + 1}b^2$ is a positive constant. Multiple (4.21) by 
$f^n$ we
obtain
$$
f^n {{\partial^2 u}\over {\partial t^2}} + nf^{n - 1} {{\partial f}\over
{\partial t}} {{\partial u}\over {\partial t}} + <\bigtriangledown_g f^{n -
2} \,,
\bigtriangledown_g u >_g + f^{n - 2} \Delta_g u + c^2 f^n u \le 0\,.\leqno
(4.22)
$$
Fix a $t > t_o$ and integrate (4.22) with respect to the \Rm $g$ and apply
Green's identity, we have
$$
\int_N f^n {{\partial^2 u}\over {\partial t^2}} dv_g + \int_N nf^{n - 1}
{{\partial f}\over {\partial t}} {{\partial u}\over {\partial t}} dv_g + c^2
\int_N f^n u dv_g \le 0
$$
for all $t > t_o$, or
$$
{{d}\over {d t}} (\int_N f^n {{\partial u}\over {\partial t}}
dv_g) \le -c^2 {\cal F} (t)\,, \leqno (4.23)
$$
where 
$$
{\cal F} (t) = \int_N f^n u dv_g > 0 \ \ \ \ {\mbox{for}} \ \ \ t > t_o\,.
$$
We have
\begin{eqnarray*}
{{d {\cal F}}\over {dt}} & = & \int_N n f^{n - 1} {{\partial f}\over
{\partial t}} u dv_g + \int_N f^n {{\partial u}\over {\partial t}} dv_g\,,\\
{{d^2 {\cal F}}\over {d t^2}} & = & \int_N n (n - 1) f^{n - 2}
|{{\partial f}\over {\partial t}}|^2 u dv_g + \int_N n f^{n - 1} 
{{\partial^2 
f}\over {\partial t^2}} u dv_g + \int_N n f^{n - 1} {{\partial  f}\over 
{\partial
t}}{{\partial  u}\over {\partial t}} dv_g\\
& \ & \ \  + {{d}\over {d t}} (\int_N f^n
{{\partial u}\over {\partial t}} dv_g)\,.
\end{eqnarray*}
Using (4.19), (4.20) and (4.23) we have
$$
{{d^2 {\cal F}}\over {d t^2}} (t) \le {{n (n - 1) C_1^2}\over {t^2}} 
{\cal F}
(t)  + {{n C_2}\over {t^2}} {\cal F} (t) + {{n C C_1}\over {t}} {\cal F} 
(t) 
- c^2 {\cal F} (t)
\leqno (4.24)
$$
for $t > t_o$. Thus we can find positive constants $\overline t > t_o$ 
and $c'$
such that 
$$
{\cal F}'' (t) \le - c'^2 {\cal F} (t) \ \ \ \ {\mbox{for \ all}} \ \ \ t 
> \overline t\,.
\leqno (4.25)
$$ 
As in the proof of theorem 4.8, (4.25) implies that ${\cal F} (t) \le 0$ 
when $t$ is
large enough.\qed 
\hspace*{0.5in}It follows as in theorem 4.8 that we can relax the 
condition on
the scalar curvature of $\og$ to
$$R_{\og} \le - {{b^2}\over {t^\alpha}}$$
for all $t > t'$ and $x \in N$, where $t' > a$, $b$ and $\alpha < 1$ are 
positive
constants. Then (4.24) gives
$$ {{d^2 {\cal F}}\over {d t^2}} (t) \le {{n (n - 1) C_1^2}\over {t^2}} 
{\cal F}
(t)  + {{n C_2}\over {t^2}} {\cal F} (t) + {{n C C_1}\over {t}} {\cal F} 
(t) - {{c^2}\over
{t^\alpha}} {\cal F} (t)\,.
$$
For $t$ large enough, since $\alpha < 1$, we have
$${\cal F}'' (t) \le - {{C'}\over {t^\alpha}} {\cal F} (t) 
\le - {{c'}\over 4} {1\over {t^2}} {\cal F}
(t)\,,$$
where $C'$ and $c' > 1$ are positive constants. The proof of lemma 1.8 
implies
that ${\cal F} (t) = 0$ for some large $t$. 

\vspace{0.3in}

{\bf \Large Appendix}

\vspace{0.3in}

{\bf Proposition A.1.} \ \ {\it For $n \ge 2$, let $\overline M$ be a 
compact $(n
+ 1)$-manifold with boundary $\partial M$ and interior $M$. Given any
smooth function $R
\in C^\infty (\overline M)\,,$ there is a \Rm $g$ (non-complete if 
$\partial M
\not= \emptyset$) defined on
$\oM$ such that $R$ is the scalar curvature of $g$ in $M$.}\\[0.1in]
{\bf Proof.} \ \ If $\oM$ is a compact manifold with boundary, then the
double of
$\oM$, defined by
$$2 \oM = (\oM \times \{1\} \cup \oM \times \{2\}\,)/\partial M\,,$$
can be given a $C^\infty$ structure as a compact manifold without 
boundary, such
that the inclusions
\begin{eqnarray*}
i_k : \oM & \hookrightarrow &2 \oM\\ 
x & \to & (x, k)/\sim\,,  \ \ \ \ \ k = 1, 2\,,
\end{eqnarray*}
are diffeomorphisms onto their range [9]. By Seeley's extension
theorem [9], any smooth function $R$ defined on $\oM \cong i_1 (\oM)$ can be
extended to a smooth function $R'$ on $2\oM$. We can modify the function $R'$
on $i_2 (M)$ so that it is negative somewhere there. Then the classification
theorem of Kazdan and Warner implies that there is a \Rm $g'$ on $2\oM$ such
that the scalar curvature of $g'$ is the extended function $R'$. We can take
$g$ to be the restriction of $g'$ on $i_1 (M)$. The \Rm $g$ is not complete
if $\partial M$ is nonempty.\qed

\pagebreak

\centerline{\bf \LARGE References}

\vspace{0.3in}

\ [1] P. Aviles \& R. McOwen, {\it Conformal deformation to constant negative
scalar curvature on noncompact Riemannian manifolds,} J. Diff. Geom. 
{\bf 27} (1988), 225-239.\\

\ [2] A. Besse, {\it Einstein Manifolds,} Springer-Verlag, Berlin-New 
York, 
1987.\\

\ [3] J. Bland \& M. Kalka, {\it Negative scalar curvature metrics on 
non-compact
manifolds,} Trans. A.M.S. {\bf 316} (1989), 433-446.\\

\ [4] F. Dobarro \& E. Lami Dozo, {\it Scalar curvature and warped 
products of
Riemannian manifolds,} Trans. A.M.S. {\bf 303} (1987), 161-168.\\

\ [5] M. Gromov \& H. B. Lawson, {\it Positive scalar curvature and the Dirac
operator on complete Riemannian manifolds,} Publ. Math. I.H.E.S. {\bf 58}
(1983), 295-408.\\

\ [6] H.B. Lawson \& M. Michelsohn, {\it Spin Geometry,} Princeton University
Press, Princeton, 1989.\\

\ [7] M.C. Leung, {\it Conformal scalar curvature equations on complete
manifolds,} Comm. in P.D.E. {\bf 20} (1995), 367-417.\\ 

\ [8] Y. Li \& W. Ni, {\it On conformal scalar curvature equations in ${\bf
R}^n$,} Duke Math. J. {\bf 57} (1988), 895-924.\\

\ [9] R.B. Melrose, {\it Differential Analysis on Manifolds with 
Corners,} M.I.T.
lecture notes, 1988.\\ 

[10] A. Ratto, M. Rigoli \& L. V\'eron, {\it Scalar curvature and conformal
deformation of hyperbolic space,} J. Functional Analysis {\bf 121} (1994),
15-77.\\

\end{document}